%% file: infinity_boundary_regularity_polynomial_degenerate_pde.tex
\documentclass[a4paper, 11pt, reqno]{amsart}

\input{header}

\title{\texorpdfstring{The Wiener Criterion at $\infty$ for Degenerate Elliptic Equations}{The Wiener Criterion at INFINITY for Degenerate Elliptic Equations}}

\author[U. G. Abdulla]{Ugur G. Abdulla}
\address[U. G. Abdulla]{Okinawa Institute of Science and Technology, Analysis and Partial Differential Equations Unit, 1919-1 Tancha, Onna-son, 904-0495, Okinawa, Japan.}
\email{ugur.abdulla@oist.jp}

\author[D. Brazke] {Denis Brazke}
\address[D. Brazke]{Okinawa Institute of Science and Technology, Analysis and Partial Differential Equations Unit, 1919-1 Tancha, Onna-son, 904-0495, Okinawa, Japan.}
\email{denis.brazke@oist.jp}

\begin{document}

	\begin{abstract}
		This paper establishes a Wiener criterion at $\infty$ to characterise the unique solvability of the Dirichlet problem for degenerate elliptic equations with power-like weights in arbitrary open sets. In the measure-theoretical context, the criterion determines whether the $\A$-harmonic measure of $\infty$ is null or positive. From the topological point of view, it presents a test for the thinness of the exterior set at $\infty$ in the $\A$-fine topology. 
	\end{abstract}

\maketitle

\tableofcontents

\section{Introduction and Main Result}

    In this paper, we prove a criterion for the unique solvability of the Dirichlet problem in arbitrary open sets $\Omega\subset \IRn$ for the degenerate elliptic equation
        \begin{align} \label{eq:dirichlet_problem} \tag{DP}
            \left\{ \begin{array}{rll}
                \A u \coloneqq \nabla \cdot (A \nabla u)    & = 0   & \quad \text{in } \Omega,  \\[6pt]
                 u                                          & = f   & \quad \text{on } \partial \Omega,
            \end{array} \right.
        \end{align}
    where $f \colon \partial \Omega \longrightarrow \IR$ is a bounded Borel measurable function, and the coefficient matrix $A \colon \IR^n \longrightarrow \IR^{n \times n}$ is measurable, symmetric and satisfies 
        \begin{align} \label{eq:ellipticity_condition}
            \frac 1\lambda \, |x|^\gamma \, |\xi|^2 \leq \xi \cdot A(x) \xi \leq \lambda \, |x|^\gamma \, |\xi|^2  && \text{for all } x, \xi \in \IRn,
        \end{align}
    for some $\lambda\geq 1$, $\gamma > 2 - n$. We characterise the existence and uniqueness of a bounded solution to \eqref{eq:dirichlet_problem} in terms of the Wiener criterion for the regularity of the boundary point at infinity ($\infty$). On one hand, it is a generalisation of the Wiener test for the regularity of $\infty$ for harmonic functions \cite{Abdulla:2007, Abdulla:2012}. On the other hand, it is the counterpart at $\infty$ of the Wiener test for the boundary regularity of solutions to degenerate elliptic equations \cite{FaJeKe:1982}.

    \medskip
    
    In order to formulate our main result, we first introduce some terminology. Consider the one point compactification $\overline \IRn \coloneqq \IRn \cup \{\infty\}$. Let $\Omega \subset \IRn$ be open and denote by $\pmb \Omega \subset \overline \IRn$ the compactification of $\Omega$. We make the convention that $\infty \in \partial \pmb \Omega$ for all unbounded open sets $\Omega \subset \IRn$. The Radon measure canonically associated with the weight $|\filll|^\gamma$ is denoted by $\mu_\gamma$:
        \begin{align} \label{eq:defi_mu_gamma}
            \mu_\gamma(E) \coloneqq \int_E |x|^\gamma \dd x,
        \end{align}
    i.e. $\mathrm d\mu = |x|^\gamma \dd x$, where $\mathrm dx$ denotes the $n$-dimensional Lebesgue measure. We say that $u \in W^{1,2}_\loc(\Omega,\mu_\gamma)$ is $\A$-harmonic if $\mathcal Au = 0$ in $\Omega$ in the weak sense, i.e.
        \begin{align}
            \int_\Omega \nabla \varphi(x) \cdot A(x) \nabla u(x) \dd x = 0 && \text{for all } \varphi \in C_c^\infty(\Omega).
        \end{align}
    We say that $v \colon \Omega \longrightarrow \IR \cup\{+\infty\}$ is $\A$-superharmonic, if
        \begin{itemize}
            \item $v$ is finite valued in a dense set of $\Omega$.
            \item $v$ is lower semi-continuous.
            \item For every open and bounded set $U \subset \Omega$ and for every $\A$-harmonic function $u \in W_\loc^{1,2}(U,\mu_\gamma) \cap C^0(\overline U)$, the inequality $v \geq h$ on $\partial U$ implies $v \geq h$ in $U$.
        \end{itemize}
    We say that $w \colon \Omega \longrightarrow \IR \cup\{-\infty\}$ is $\A$-subharmonic if $-w$ is $\A$-superharmonic. Throughout this article, by a solution of \eqref{eq:dirichlet_problem} we mean a Perron–Wiener–Brelot solution (PWB solution for short). Assuming that boundary function $f \in C^0(\partial \pmb\Omega)$, the upper PWB-solution $\overline \H_f^\Omega$ (resp. lower PWB-solution $\underline \H_f^\Omega$) is defined as 
        \begin{align}
            \overline \H^\Omega_f(x) \coloneqq \inf \{v(x)\} \quad (\text{resp. } \underline \H^\Omega_g(x) \coloneqq \sup \{w(x)\}),   && \text{for all } x \in \Omega,
        \end{align}
    where the infimum (resp. supremum) is taken over all $\A$-superharmonic functions (resp. $\A$-subharmonic functions) such that 
        \begin{align}\label{eq:defi_bc}
            \liminf_{\substack{y \to x \\ y \in \Omega}} v(y) \geq f(x) \quad \big(\text{resp. } \limsup_{\substack{y \to x \\ y \in \Omega}} v(y) \leq f(x) \big) && \text{for all } x \in \partial \pmb \Omega
        \end{align}
    We say that a boundary function $f \colon \partial \pmb\Omega \longrightarrow \IR$ is $\A$-resolutive, if
        \begin{align}
            \underline \H^\Omega_f \equiv \overline \H^\Omega_f \eqqcolon \H_f^\Omega,
        \end{align}
    and the latter is called a generalised solution (or PWB solution) of \eqref{eq:dirichlet_problem}. Perron's method and its Wiener and Brelot refinements imply that continuous boundary functions, and therefore bounded Borel measurable functions, form a subclass of resolutive boundary functions \cite{Doob:1984} (see Section \ref{sec:preliminaries}).

    \medskip
    
    A generalised solution is $\A$-harmonic in $\Omega$ and unique by construction. However, it is accomplished by prescribing {\it the behaviour of the solution at $\infty$}, i.e. by enforcing \eqref{eq:defi_bc} at the boundary point $\infty$. The theory, while identifying a class of unique solvability, leaves the following question open: {\it would a unique bounded solution exist without prescribing its behaviour at $\infty$?} The principal goal of this paper is to answer this question. {\it Furthermore, without loss of generality, we assume that the open set $\Omega$ has at least one unbounded connected component.} It is evident that, without this assumption, within each bounded component, the Dirichlet problem is uniquely solvable. This guarantees the unique solvability of \eqref{eq:dirichlet_problem} in the entire open set $\Omega$ without prescribing the solution at $\infty$. In fact, for any open set $\Omega$ satisfying this condition, without prescribing the behaviour of the solution at $\infty$, there exists one and only one, or infinitely many solutions (see Proposition \ref{prop:harmonic_measure::uniqueness}).

    \medskip
    
    By fixing any real number $\overline{f} \in \IR$, we can assign $f(\infty) = \overline{f}$ and construct a unique generalised solution $\H_f^\Omega$. The major problem is now the following: {\it How many bounded solutions do we actually have, or does the constructed solution depend on $\overline{f}$?.} 

    \medskip
    
    To formulate the corresponding measure-theoretical problem, we recall the concept of $\A$-harmonic measure of Borel subsets of $\partial \pmb \Omega$. For a given Borel set $B \subset \partial\pmb\Omega$, denote as $\mathbf 1_B$ the indicator function of $B$. The $\A$-harmonic measure of $\infty$ is defined as
        \begin{align}
            \mu_{\pmb \Omega}(\filll,\infty)=\H_{\mathbf 1_\infty}^\Omega(\filll).
        \end{align}
    We say that $\infty$ is an $\A$-harmonic measure null set, if $\mu_{\pmb\Omega}(\filll,\infty) \equiv 0$ in $\Omega$. Otherwise, $\infty$ is a set of positive $\A$-harmonic measure. The measure-theoretical counterpart of the major problem is to characterise geometrically whether $\infty$ {\it is of null or positive $\A$-harmonic measure}.
    
    To characterise the uniqueness, we employ the concept of regularity of $\infty$ introduced in \cite{Abdulla:2007, Abdulla:2008, Abdulla:2012}. 

    \begin{defi} \label{defi:regularityofinfinity}
        $\infty$ is said to be regular (resp. irregular) for $\Omega$ if it is an $\A$-harmonic measure null (resp. positive) set.
    \end{defi}

    The concept of regularity of $\infty$ is fundamentally different from the classical definition of regularity of finite boundary points. Nevertheless, it provides the key insight to characterise the regularity of solutions at $\infty$. In Proposition \ref{prop:uniqueness::regular_point} below, we prove that an equivalent definition of the regularity of $\infty$ may be given depending on whether the following condition is satisfied:
        \begin{align} \label{eq:defi:regularity_limit}
           \liminf_{\substack{x \to \infty \\ x \in \partial \Omega}} f(x) \leq \liminf_{\substack{x \to \infty \\ x \in \Omega}} \H_f^\Omega(x) \leq \limsup_{\substack{x \to \infty \\ x \in \Omega}} \H_f^\Omega(x) \leq \limsup_{\substack{x \to \infty \\ x \in \partial \Omega}} f(x).
        \end{align}
    for every bounded Borel measurable boundary function $f$ and for every solution $\H_f^\Omega$ with arbitrary prescribed boundary condition $f(\infty) = \overline{f}$. Note that if $f|_{\partial\Omega}$ has a limit at $\infty$, \eqref{eq:defi:regularity_limit} means that the solution $\H_f^\Omega$ is continuous at $\infty$ by taking the limit value of $f|_{\partial\Omega}$ without being specified.

    \medskip

     To formulate the corresponding topological problem, recall that the $\A$-fine topology is defined as the coarsest topology such that every $\A$-superharmonic function in $\IRn$ is continuous. Given $E \subset \IRn$, we define $\overline E^\A$ as the closure of $E$ in the $\A$-fine topology.

    \begin{defi} \cite{Abdulla:2012}
        Let $E \subset \IRn$. We say that $E$ is $\A$-thin at $\infty$ if $E$ satisfies one of the following two conditions:
            \begin{enumerate}
                \item $E$ is bounded.
                \item $E$ is unbounded and there exists an $\A$-superharmonic function $u \colon \IRn \longrightarrow \IR$ such that 
                    \begin{align}
                        \liminf_{\substack{x \to \infty \\ x \in \IRn}} u(x) < \liminf_{\substack{x \to \infty \\ x \in E}} u(x)
                    \end{align}
            \end{enumerate}
        If $E$ is not $\A$-thin at $\infty$, we say that $E$ is $\A$-thick at $\infty$.
    \end{defi}
       
    Our principal problem is to derive a "geometric" characterisation of the regularity of $\infty$ and the uniqueness of the solution of the Dirichlet problem \eqref{eq:dirichlet_problem}. The precise characterisation depends on the sparseness of $\Omega$ at $\infty$ or, equivalently, on the fine topological thinness of the exterior set at $\infty$. To measure the critical thinness, we recall the concept of $\A$-capacity: let $B_R \coloneqq \{x \in \IRn : |x| < R\}$ be a ball and $K \subset B_R$ be compact. The $\A$-capacity of $K$ relative to $B_R$ is defined as
    
    \begin{align} \label{eq:defi_capacity}
        \capac(K,B_R) \coloneqq \inf \bigg\{ \int_{B_R} \nabla \varphi(x) \cdot A(x) \nabla \varphi(x) \dd x : \varphi \in C_c^\infty(B_R), \ \varphi \geq 1 \text{ in } K \bigg\}.
    \end{align}
    
    Our main theorem reads:
    
    \begin{thrm} \label{maintheorem} 
        The following conditions are equivalent:
            \begin{enumerate}
                \item $\infty$ is regular (resp. irregular) for $\Omega$.
                \item There exist(s) a unique (resp. infinitely many) bounded solution(s) of \eqref{eq:dirichlet_problem}.
                \item $\Omega^c$ is $\A$-thick (resp. $\A$-thin) at $\infty$. 
                \item The Wiener integral 
                    \begin{align}\label{Wienertest}
                        \int_1^\infty \frac{\capac(\Omega^c \cap \overline B_t, B_{2t})}{\capac(B_t, B_{2t})} \frac{\mathrm dt}{t}
                    \end{align}
                    diverges (resp. converges).
                \item The Wiener sum
                        \begin{align}
                            \sum_{k = 1}^\infty (2^{-k})^{n + \gamma - 2} \, \capac(E_k, B_{2^{k + 1}})
                        \end{align}
                    is divergent (resp. convergent), where $E_k \coloneqq \Omega^c \cap \{2^{k - 1} \leq |x| \leq 2^k\}$.
            \end{enumerate}
    \end{thrm}  

    We note that the Wiener test given in \textit{(5)} coincides with the Wiener test in \cite{Abdulla:2007, Abdulla:2012} in the case $\gamma = 0$ and $n \geq 3$.

    \subsection{Historical remarks} It is appropriate to make some remarks concerning the well-posedness of the Dirichlet problem for elliptic PDEs in arbitrary open sets. The strategy for solving the classical Dirichlet problem for harmonic functions in arbitrary bounded open sets may well be expressed by the citation from Wiener's celebrated paper. As pointed out by Lebesgue and independently by Wiener \cite{Wiener:1924a}, {\it "the Dirichlet problem divides itself into two parts, the first of which is the determination of harmonic functions corresponding to certain boundary condition, while the second is the investigation of the behaviour of this function in the neighbourhood of the boundary".} The existence of a unique generalised solution of the Dirichlet problem in arbitrary bounded open sets $\Omega$ with prescribed boundary values on $\partial \Omega$ is realized within the class of {\it resolutive boundary functions}, identified by Perron's method and its Wiener \cite{Wiener:1924, Wiener:1924a} and Brelot \cite{Brelot:1971} refinements, referred to as the PWB method. Continuous boundary functions are a subclass of resolutive boundary functions, and a finite boundary point is {\it regular} if the corresponding PWB solution takes the boundary value continuously.
    
    \medskip
    
    The regularity of a finite boundary point $x_o \in \partial\Omega$ is a problem of local nature as it depends on the measure-geometric properties of the boundary in the neighbourhood of $x_o$. It is a standard textbook example to demonstrate that an isolated boundary point is always irregular. On the other hand, any simply connected bounded open plane region is a regular set, i.e. all the boundary points are regular. Lebesgue constructed an example of a simply connected domain in $\IR^3$, which is now called a Lebesgue spine, with an irregular boundary point on a sufficiently sharp edge protruding into the interior of the region \cite{Lebesque:1912}. Wiener, in his seminal works \cite{Wiener:1924, Wiener:1924a}, proved a necessary and sufficient condition for the finite boundary point $x_o \in \partial\Omega$ to be regular in terms of the "thinness" of the complementary set in the neighbourhood of $x_o$. A key advance made in Wiener's work was the introduction of the concept of capacity, a sub-additive set function dictated by the Laplace operator. It accurately measures the thinness of the complementary set in the neighbourhood of $x_o$, and determines the regularity of the boundary point for harmonic functions. Formalised through the powerful Choquet capacitability theorem \cite{Choquet:1954}, this concept became a standard tool in potential theory to characterise boundary regularity. The Wiener criterion for the boundary continuity of harmonic functions became a canonical result, driving the boundary regularity theory for elliptic and parabolic PDEs. In \cite{LiStWe:1963}, it is proved that the Wiener criterion for the regularity of finite boundary points concerning a second-order divergence form uniformly elliptic operator with bounded measurable coefficients coincides with the classical Wiener criterion for the boundary regularity of harmonic functions. The Wiener criterion for the regularity of finite boundary points for the linear degenerate elliptic equations is proved in \cite{FaJeKe:1982}. Wiener criterion for the regularity of finite boundary points for quasilinear elliptic equations was settled due to \cite{Mazja:1976, GaZie:1977, LinMar:1985, KilMal:1994}.
    
    \medskip
    
    To solve the Dirichlet problem in an unbounded open set, Brelot introduced the idea of compactifying $\IRn$ to $\overline \IRn \coloneqq \IRn \cup \{\infty\}$, where $\infty$ is the point at infinity of $\IRn$ \cite{Bre:1944}. The PWB method is extended to the compactified framework, and provides a powerful existence and uniqueness result for \eqref{eq:dirichlet_problem} in arbitrary open sets in the class of resolutive boundary functions. The new concept of regularity of $\infty$ was introduced in \cite{Abdulla:2007} for the classical Dirichlet problem, and in \cite{Abdulla:2008} for its parabolic counterpart. The Dirichlet problem with bounded Borel measurable boundary functions has one, and only one, or infinitely many solutions without prescribing the boundary value at $\infty$. The point at $\infty$ is called regular if there is a unique solution and irregular otherwise. Equivalently, in the measure-theoretical context, the new concept of regularity resp. irregularity of $\infty$ is introduced according to whether the harmonic measure of $\infty$ is null resp. positive. In \cite{Abdulla:2007}, the Wiener criterion for the regularity of $\infty$ for the classical Dirichlet problem in an open set $\Omega \subset \IRn$ with $N\geq 3$ is proved. In \cite{Abdulla:2012}, it is proved that the Wiener criterion at $\infty$ for the linear second-order divergence form uniformly elliptic PDEs with bounded measurable coefficients coincides with the Wiener criterion at $\infty$ for the Laplace operator. The Wiener criterion at $\infty$ for the heat equation was proved in \cite{Abdulla:2008}. The main goal of this paper is to extend the results of \cite{Abdulla:2007, Abdulla:2012} to the case of degenerate elliptic PDEs with power-like weights $|x|^\gamma$, $\gamma > 2 - n$. On one hand, Theorem~\ref{maintheorem} is in direct analogy to the Wiener criterion at $\infty$ in the unweighted case \cite{Abdulla:2012}. On the other hand, Theorem~\ref{maintheorem} is the counterpart at $\infty$ of the Wiener criterion for the regularity of finite boundary points for weighted elliptic PDEs with weight $w$ being either in the Muckenhoupt class $A_2$, or $w = |\det(\mathrm D f)|^{1 - \frac{2}{n}}$, where $f: \IRn \longrightarrow \IRn$ is a global quasi-conformal map \cite{FaJeKe:1982}. This includes power-like weights $|x|^\gamma$, $\gamma > -n$. In \cite{FaJeKe:1982} it is proved that a boundary point at the origin $0 \in \partial \Omega$ is regular if and only if
        \begin{align}\label{eq:wcforfinitebp}
            \int_{0^+} \frac{1}{\capac(B_t,B_{2t})} \frac{\mathrm dt}{t} < +\infty \qquad \textbf{or} \qquad \int_{0^+} \frac{\capac(\Omega^c \cap \overline B_t,B_{2t})}{\capac(B_t,B_{2t})} \frac{\mathrm dt}{t} = +\infty.   
        \end{align}
    In the special case $w(x) = |x|^\gamma$, $\gamma > -n$, it implies that $0\in\partial\Omega$ is always regular if $\gamma < 2 - n$, but if $\gamma \geq 2 - n $ then it is regular if and only if the second condition of \eqref{eq:wcforfinitebp} is satisfied. The Wiener criterion at $\infty$ expressed in \eqref{Wienertest} is the analogue of the Wiener criterion at the finite boundary point expressed in the second integral in \eqref{eq:wcforfinitebp}.
    
    \medskip
    
    However, Theorem~\ref{maintheorem} leaves an open problem on the optimality of the restriction $\gamma > 2 - n$ for the validity of the Wiener criterion for the regularity of $\infty$. From \cite{Abdulla:2018} it follows that in the special case $n = 2$, the restriction $\gamma > 0$ is optimal. The limit case $\gamma = 0$ corresponds to the case of \eqref{eq:dirichlet_problem} for the second-order divergence form uniformly elliptic PDE with bounded measurable coefficients. In this case, $\infty$ is always regular for any unbounded Greenian open set of $\IR^2$, that is to say the Dirichlet problem always has a unique bounded solution without specifying the behaviour at $\infty$. Remarkably, in this case, the Wiener criterion characterises the removability of the fundamental singularity. More precisely, let $\Omega \subset \IR^2$ be a Greenian open set, and let $x_o$ be a boundary point (finite or $\infty$). Consider the singular Dirichlet problem for a linear divergence form uniformly elliptic operator with bounded measurable coefficients in the class $O(\log|x - x_o|)$ if $x_o$ is finite, and in the class of functions with logarithmic growth if $x_o = \infty$. In \cite{Abdulla:2018}, it is proved that the Wiener criterion at $x_o$ is a necessary and sufficient condition for the unique solvability of the singular Dirichlet problem, and equivalently for the removability of the logarithmic singularity. Precisely, in \cite{Abdulla:2018}, the concept of log-regularity resp. log-irregularity of the boundary point (finite or $\infty$) is introduced according to if its log-harmonic measure is null resp. positive, and the removability of the logarithmic singularity is expressed in terms of the Wiener criterion for the log-regularity of $x_o$. We address the problem of characterising the regularity of $\infty$ for arbitrary unbounded open sets in the \eqref{eq:dirichlet_problem} for the weight parameter $\gamma$ in the range $-n < \gamma \leq 2 - n$, as well as the problem of characterising singularities for degenerate elliptic PDEs in a subsequent work. 
    
    \subsection{Notation} We mainly use standard notation in this article. We denote by $\IN$ the set of positive integers, and for $n \in \IN$, $n \geq 2$, we denote the Euclidean space by $\IRn$. The complement of a set $\Omega \subset \IRn$ is denoted by $\Omega^c \coloneqq \IRn \setminus \Omega$, and the closure is denoted by $\overline \Omega$. Given $x,y \in \IRn$, we denote the standard scalar product by $x \cdot y$ and define $|x|^2 \coloneqq x \cdot x$. We denote $\operatorname{diag} \coloneqq \{(x,x) \in \IRn \times \IRn : x \in \IRn\}$. If there is no confusion with the underlying domain and limit, we denote by $f^*$ resp. $f_*$ the limsup resp. liminf of the function $f$ at $\infty$. We write $X \lesssim Y$ if there exists a constant $c = c(n,\gamma) > 0$ such that $X \leq cY$, where $n$ is the dimension of the underlying space and $\gamma$ is the exponent for the weight in the ellipticity condition \eqref{eq:ellipticity_condition}. We write $X \approx Y$ if $X \lesssim Y$ and $Y \lesssim X$.

    \medskip
    
    \textbf{Acknowledgements:} D.B. thanks Daniel Tietz and Chenming Zhen for highly valuable discussions.

\section{Preliminaries} \label{sec:preliminaries}

    Given $\gamma > -n$, we define the measure $\mu_\gamma$ as in \eqref{eq:defi_mu_gamma}. Since the function $x \longmapsto |x|^\gamma$ is $2$-admissible (in the sense of \cite[Chapter 1]{HeKiMa:2006}), the measure $\mu_\gamma$ is doubling and $\mu_\gamma$ and the Lebesgue measure are mutually absolutely continuous (see \cite[Example 1.6 and Chapter 15]{HeKiMa:2006}). A direct computation shows that
        \begin{align} \label{eq:measure_of_ball}
            \mu_\gamma(B_r(0)) = \sigma_{n,\gamma} r^{n + \gamma}   && \text{for all } r > 0, 
        \end{align}
    where $\sigma_{n,\gamma} \coloneqq \frac{n}{n + \gamma} \mu_0(B_1(0))$. From this computation we obtain for $x \in \IRn$ and $r > |x|$ the following rough estimate
        \begin{align} \label{eq:estimate_measure_of_ball}
            (r - |x|)^{n + \gamma} \lesssim \mu_\gamma(B_r(x)) \lesssim (r + |x|)^{n + \gamma}.
        \end{align}
    Recall that $A \colon \IRn \longrightarrow \IR^{n \times n}$ is a measurable and symmetric coefficient matrix that satisfies the ellipticity condition \eqref{eq:ellipticity_condition} for some $\gamma > 2 - n$. We remark that the constant $\gamma$ is not covering the full range of parameters that make the weight $\mu_\gamma$ 2-admissible (in the sense of \cite[Chapter 1]{HeKiMa:2006}). The main reason for this restriction is the existence and decay of full space potentials of the weighted equation $\A u \coloneqq \div(A \nabla u) = 0$ (see Section \ref{sec:full_space_potential_theory}), which can only be guaranteed if $\gamma > 2 - n$. 

    \medskip
    
    We define the upper resp. lower PWB solution as in the introduction. Both functions are $\A$-harmonic (see \cite[Theorem 9.2]{HeKiMa:2006}). Using the comparison principle, we find that $\underline \H^\Omega_g \leq \overline \H^\Omega_g$ for every boundary function $g$. Moreover, both functions are order-preserving, i.e.
        \begin{align}
            \underline \H^\Omega_g \leq \underline \H^\Omega_{\tilde g} && \overline \H^\Omega_g \leq \overline \H^\Omega_{\tilde g}
        \end{align}
    for every bounded Borel measurable functions $g, \tilde g \colon \partial \pmb \Omega \longrightarrow \IR$ such that $g \leq \tilde g$.

    \medskip
    
    We say that a boundary function $g \colon \partial \pmb\Omega \longrightarrow \IR$ is $\A$-resolutive, if $\underline \H^\Omega_g = \overline \H^\Omega_g$ and are not constant $\pm \infty$. It is well known that every continuous boundary function, and therefore also every bounded Borel measurable function $g \colon \partial \pmb \Omega \longrightarrow \IR$ is $\A$-resolutive due to the linearity of the equation. For the reader's convenience, we provide a sketch of the proof of this fact:

    \begin{lemma} 
        Let $\Omega \subset \IRn$ be an $\A$-Greenian open set, i.e. it possesses a non-trivial and non-negative $\A$-superharmonic function. Then every bounded and Borel measurable function $f \colon \partial \pmb\Omega \longrightarrow \IR$ is $\A$-resolutive.
    \end{lemma}
        \proofb The PWB solution is stable under monotone convergence, see \cite[Section VIII.6 (e)]{Doob:1984}. The result is only stated for the Laplace equation, but it is also true for weighted elliptic equations, as the argument carries over to that case almost verbatim. Using the resolutivity of continuous functions \cite[Theorem 9.25]{HeKiMa:2006} and aforementioned argument of \cite[Section VIII.6 (e)]{Doob:1984}, we find that indicator functions of Borel measurable subsets of $\partial \pmb \Omega$ are resolutive. Using the fact that the space of $\A$-resolutive functions is a vector space, including \cite[Section VIII.6 (e)]{Doob:1984} indicator functions of Borel measurable boundary subsets, the claim follows. \proofe

    \subsection{Preliminaries in Potential Theory} Let $B \subset \IRn$ be a ball and $K \subset B$ be compact. We define $\capac(K,B)$ as in \eqref{eq:defi_capacity}. For arbitrary sets $E \subset B$, we define
        \begin{align}
            \capac(E,B) \coloneqq \inf_{\substack{E \subset U \subset B \\ U \text{ open}}} \sup_{\substack{K \subset U \\ K \text{ compact}}} \capac(K,B).
        \end{align}
    The set function $E \longmapsto \capac(E,B)$ is a Choquet capacity (see \cite[Theorem 2.5]{HeKiMa:2006}). We say that a property holds quasi-everywhere if it holds up to a set of capacity $0$ (see \cite[Section 2.7]{HeKiMa:2006}). We mention that the set of irregular finite boundary points of any open set is a set of capacity $0$ (see \cite[Theorem 8.10]{HeKiMa:2006})

    \begin{lemma} \label{lemma:capacity_properties} \text{}
        \begin{enumerate}
            \item $\capac(E,B) \leq \capac(E',B)$ for all $E \subset E' \subset B$.
            \item $\capac(E,B) \leq \capac(E,B')$ for all $E \subset B' \subset B$.
            \item $\capac(B_r(x), B_{2r}(x)) = \capac(\overline B_r(x), B_{2r}(x))$.
            \item \label{it:lemma:capacity:properties:measure_to_capacity} $\mu_\gamma(B_r(x)) \approx r^2 \capac(B_r(x), B_{2r}(x))$ for all $x \in \IRn$, $r > 0$.
            \item $\capac(E, B_{2t}(x)) \approx \capac(E,B_{2s}(x))$ for all $ 0 < t < s \leq 2t$, $E \subset B_t(x)$.
        \end{enumerate}
    \end{lemma}

    The proofs can be found in \cite[Chapter 2]{HeKiMa:2006}.
        
    \medskip
           
    It is well known that there exists a unique minimiser $\hat \R^1[K,B] \in W_0^{1,2}(B,\mu_\gamma)$ of the functional minimised in \eqref{eq:defi_capacity}, which we call the $\A$-capacitary potential of $K$ in $B$. Moreover, $0 \leq \hat \R^1[K,B] \leq 1$ almost everywhere, $\hat \R^1[K,B]$ is $\A$-harmonic in $B \setminus K$ and $\hat \R^1[K,B]$ is $\A$-superharmonic in $B$. Without loss of generality we always choose the $L^1$-representative of $\hat \R^1[K,B]$ such that
        \begin{align}
            \hat \R^1[K,B](x) = \Essliminf_{y \longrightarrow x} \hat \R^1[K,B](y).
        \end{align}
    We note that $\hat \R^1[K,B]|_{B \setminus K} = \H^{B \setminus K}_f$ where $f \coloneqq \varphi|_{\partial (B\setminus K)}$, and where $\varphi \in C_c^\infty(B)$ such that $\varphi = 1$ in an open neighbourhood of $K$ in $B$ (see \cite[Corollary 9.29]{HeKiMa:2006}). For more details, see \cite[Section 6.16]{HeKiMa:2006} (see also \cite[Theorem 1.20]{FaJeKe:1982}).

    \medskip

    Let $\nu$ be a Radon measure on $B$. We say that $u \in L^1(B,\mu_\gamma)$ is a weak solution to $\A u = \nu$ vanishing on $\partial B$ if
        \begin{align}
            \int_B u(x) \, \Phi(x) \, |x|^\gamma \dd x = \int_B \G[\Phi  \,|\cdot |^\gamma] \dd \nu
        \end{align}
    for all $\Phi \in L^\infty(B,\mu_\gamma)$, where $\G \colon W^{-1,2}(B,\mu_\gamma) \longrightarrow W_0^{1,2}(B,\mu_\gamma)$ is the Greens operator defined on the dual space $W^{-1,p'}(B,\mu_\gamma) \coloneqq [W_0^{1,p}(B,\mu_\gamma)]'$ with $p' = \frac{p}{p - 1}$. Note that $\Phi |\cdot|^\gamma \in W^{-1,p}(B,\mu_\gamma)$ for all $p > 1$, from which we obtain that $\G[\Phi w] \in C^0(\overline B)$. It is well known that a unique weak solution $u \in L^1(B,\mu_\gamma)$ exists. Moreover, $u \geq 0$ almost everywhere and there exists $p_0 < 2n$ such that $u \in W_0^{1,p}(B,\mu_\gamma)$ for all $1 \leq p < p_0$. In particular $u \in W_0^{1,2}(B,\mu_\gamma)$ if $\nu \in W^{-1,2}(B, \mu_\gamma)$. See \cite[Chapter 2]{FaJeKe:1982} for more details.

    \medskip

    With the notion of weak solutions, we define the Greens function $g \colon (B \times B) \setminus \operatorname{diag} \longrightarrow \IR$ as the unique weak solution of the equation $\A g(\filll,y) = \delta_y$ vanishing on $\partial B$ for all $y \in B$. We note that $g(\filll, y)$ is $\A$-harmonic in $B \setminus \{y\}$, $g \geq 0$ and $g(\filll,y) \in W^{1,2}(B \setminus B_r(y), \mu_\gamma)$ for all $r > 0$. In addition, the Greens function $g$ is symmetric (see \cite[Proposition 2.8]{FaJeKe:1982}) and is estimated via (see Lemma \ref{lemma:capacity_properties} \textit{(3)})
        \begin{align} \label{eq:greens_function_estimate_local}
            g(x,y) \approx \int_{|x - y|}^R \frac{1}{\capac(B_t(x)), B_{2t}(x))} \frac{\mathrm dt}{t}   && \text{for all } x,y \in \tfrac 14B.
        \end{align}
   where $R > 0$ is the radius of $B$ (see \cite[Theorem 3.3]{FaJeKe:1982}).
    
    \medskip
    
    Given a Radon measure $\nu$ with compact support in $B$, the function $u \colon B \longrightarrow \IR$ defined via
        \begin{align} \label{eq:greens_representation_formula}
            u(x) \coloneqq \int_B g(x,y) \dd \nu(y) && \text{for all } x \in B
        \end{align}
    is the unique weak solution of $\A u = \nu$ vanishing on $\partial B$ (see \cite[Lemma 2.7]{FaJeKe:1982}). In particular, there exists a unique $\nu \in W^{-1,2}(B,\mu_\gamma)$ such that $\A \hat \R^1[K,B] = \nu$ in $B$ and $\nu(B) = \capac(K)$ (see \cite[Proposition 1.22]{FaJeKe:1982}). Moreover $\supp(\nu) \subset \partial K$. This measure is called the capacitary measure of $K$ in $B$ and the potential $\hat \R^1[K,B]$ obeys \eqref{eq:greens_representation_formula}.

\subsection{Potential Theory in the whole space} \label{sec:full_space_potential_theory}

    The preceding section discussed potential theory for functions defined on a ball $B \subset \IRn$. In this section, we want to use these preliminaries to extend the results to the whole space setting. A simple calculation using Lemma \ref{lemma:capacity_properties} shows that
        \begin{align} \label{eq:structural_capacity_estimate}
            |x - y|^{2 - (n + \gamma)} \approx \int_{|x - y|}^\infty \frac{1}{\capac(B_t(x),B_{2t}(x))} \frac{\mathrm dt}{t}  && \text{for all } |x| > 4 |y|.
        \end{align}    
    With this estimate at hand, we can define whole space potentials. More precisely, given a compact set $K \subset \IRn$, we construct $\hat \R^1[K,\IRn] \colon \IRn \longrightarrow \IR$ as the limit of the increasing sequence of $\hat \R^1[K,B_R]$ as $R \to \infty$. Since $\hat \R^1[K,B_R] \leq 1$ for all $R > 0$, again the function $\hat \R^1[K,\IRn]$ is $\A$-harmonic in $\IRn \setminus K$ by the Harnack convergence theorem, and $\A$-superharmonic in $\IRn$. In addition, we obtain a capacitary upper bound on the potential $\hat \R^1[K,\IRn]$ in the following way. Let $K \subset B_r$, let $G_R \colon B_R \times B_R \longrightarrow \IR$ be the Greens function on $B_R$ and let $\nu_R$ be the capacitary potential of $\hat \R^1[K,B_R]$. From the Greens representation formula \eqref{eq:greens_representation_formula} we find for $x \in \IRn \setminus B_{2r}$
        \begin{align}
\label{eq:potential_capacitary_bound1}  \hat \R^1[K,\IRn](x)    &  = \lim_{R \to \infty} \int_{\partial K} G_R(x,y) \dd \nu_R(y)  \\[6pt]
                                                                & \lesssim \liminf_{R \to \infty} \int_{\partial K} \int_{|x - y|}^\infty \frac{1}{\capac(B_t(x),B_{2t}(x))} \frac{\mathrm dt}{t} \dd \nu_R(y) \\[6pt]
\label{eq:potential_capacitary_bound2}                          & \leq \capac(K,B_{2r}) \, \int_r^\infty \frac{1}{\capac(B_t(x),B_{2t}(x))} \frac{\mathrm dt}{t}.
        \end{align}
    In particular, $\hat \R^1[K,\IRn](x) \to 0$ as $x \to \infty$. Moreover, every open set of $\IRn$ is $\A$-Greenian since full space potentials provide a non-trivial and non-negative $\A$-superharmonic function.
    
    \medskip
    
    In contrast to the unweighted case $\gamma = 0$, a similar lower bound (in terms of capacity) is generally unavailable in the whole space setting due to a lack of a weighted Sobolev inequality. However, it follows from the weighted Sobolev inequality proved in \cite{CabRos:2013} that a suitable capacitary lower bound is available in the case $\gamma > 0$. As this fact is not required for our purposes, we omit the details.

    \medskip

    Similar to the unweighted case (\cite[Lemma 2.1]{Abdulla:2007}, \cite[page 3391]{Abdulla:2012}), in the following lemma we use full space potentials to show that we are always able to construct a generalised solution of \eqref{eq:dirichlet_problem} such that \eqref{eq:defi:regularity_limit} is satisfied.

    \begin{lemma} \label{lem:construction_generalised_solution_good_limit}
        Let $f \colon \partial \Omega \longrightarrow \IR$ be a bounded Borel measurable function. Let $\overline f \in \IR$ such that $f_* \leq \overline f \leq f^*$ and extend $f(\infty) = \overline f$. Then $\H^\Omega_f$ satisfies \eqref{eq:defi_bc}.
    \end{lemma}
        \proofb We denote $u \coloneqq \H^\Omega_f$. Let $M \coloneqq \|f\|_{\sup}$ and let $\varepsilon > 0$. Then $|u| \leq M$ since the PWB-method is order-preserving and we find $R > 0$ such that
            \begin{align}
                f_* - \varepsilon \leq f(y) \leq f^* + \varepsilon  && \text{for all } y \in \partial \Omega \cap B_R^c.
            \end{align}
        We now define the function $w_{\pm} \colon \Omega \cap \overline B_R^c \longrightarrow \IR$ via
            \begin{align}
                w_+ \coloneqq f^* + \varepsilon + 2M \hat \R^1[B_R,\IRn],   && w_- \coloneqq f_* - \varepsilon - 2M \hat \R^1[B_R,\IRn].
            \end{align}
        Since the PWB-method is order-preserving, we find that $w_- \leq u \leq w_+$ in $\Omega \cap B_R^c$. In the limit $x \to \infty$ and then $\varepsilon \to 0$, we find that $f_* \leq u_* \leq u^* \leq f^*$, which was the claim. \proofe

    \subsection{Wiener Integral} \label{subsec:wiener_integral} We define the Wiener integral
        \begin{align}
            \WW(\delta, \Omega) \coloneqq \int_\delta^\infty \frac{\capac(\Omega^c \cap \overline B_t(x), B_{2t}(x))}{\capac(B_t(x),B_{2t}(x))} \frac{\mathrm dt}{t}
        \end{align}
    for $\Omega \subset \IRn$ open, $x \in \IRn$ and $\delta > 0$. We say that the Wiener integral converges if there exists $x \in \IRn$ and $\delta > 0$ such that $\WW(\delta, \Omega) < \infty$. Otherwise, we say that the Wiener integral diverges. In addition, a simple computation using Lemma \ref{lemma:capacity_properties} shows that convergence resp. divergence of the Wiener integral is independent of the centre points $x \in \IRn$. Therefore, without loss of generality, we always take $x = 0$. Furthermore, again using Lemma \ref{lemma:capacity_properties}, convergence resp. divergence of the Wiener integral is equivalent to convergence resp. divergence of the Wiener integral when tested with open balls instead of closed balls. We observe that the Wiener integral is non-increasing:
        \begin{align}
            \WW(\delta,\Omega')   & \leq \WW(\delta,\Omega) && \text{for every } \delta > 0, \ \Omega \subset \Omega'\text{ open},  \\[6pt]
            \WW(\delta',\Omega)   & \leq \WW(\delta,\Omega) && \text{for every } 0 < \delta < \delta', \ \Omega \subset \IRn \text{ open}.
        \end{align}
    Moreover it follows from Lemma \ref{lemma:capacity_properties} that if there exists $\delta^* > 0$ such that $\WW(\delta^*,\Omega) < \infty$, then $\WW(\delta,\Omega) < \infty$ for all $\delta > 0$. Similarly, this holds for divergence of the Wiener integral. 

    \medskip

    We note that the Wiener integral is equivalent to the Wiener sum tested with annuli instead of balls. This version of the Wiener sum is used throughout \cite{Abdulla:2007, Abdulla:2012} (see also \cite{Abdulla:2018}) and takes on a similar role in this article. More precisely, convergence resp. divergence of the Wiener integral is equivalent to convergence resp. divergence of
        \begin{align}
            \WW^\Sigma(m,\Omega) \coloneqq \sum_{k = m}^\infty (2^{-k})^{n + \gamma - 2} \capac(\Omega^c \cap A_k,B_{2^{k + 1}}),
        \end{align}
    for some $m \in \IN$, where $A_k \coloneqq \overline B_{2^k} \setminus B_{2^{k - 1}}$. To see this, we use Lemma \ref{lemma:capacity_properties} to find that
        \begin{align}
            \WW^\Sigma(m,\Omega)    & \leq \sum_{k = m}^\infty (2^{-k})^{n + \gamma - 2} \capac(\Omega^c \cap \overline B_{2^k},B_{2^{k + 1}})  \\[6pt]
                                    & \lesssim \sum_{k = m}^\infty (2^{-k})^{n + \gamma - 2} \capac(\Omega^c \cap \overline B_{2^k},B_{2^{k + 2}}) \\[6pt]
                                    & \lesssim \sum_{k = m}^\infty (2^{-k})^{n + \gamma - 2} \capac(\Omega^c \cap \overline B_{2^k},B_{2^{k + 2}}) \int_{2^k}^{2^{k + 1}} \frac{\mathrm dt}{t}  \\[6pt]
                                    & \lesssim \sum_{k = m}^\infty \int_{2^k}^{2^{k + 1}} \frac{\capac(\Omega^c \cap \overline B_t,B_{2t})}{t^{n + \gamma - 2}} \frac{\mathrm dt}{t}    \\[6pt]
                                    & \lesssim \sum_{k = m}^\infty \int_{2^k}^{2^{k + 1}} \frac{\capac(\Omega^c \cap \overline B_t,B_{2t})}{\capac(B_t,B_{2t})} \frac{\mathrm dt}{t} = \WW(2^m,\Omega).
        \end{align}
   To show the reverse inequality, we use Lemma \ref{lemma:capacity_properties} to compute
        \begin{align}
            \WW(2^m,\Omega) & \lesssim \int_{2^m}^\infty \frac{\capac(\Omega^c \cap \overline B_t,B_{2t})}{t^{n + \gamma - 2}} \frac{\mathrm dt}{t}    \\[6pt]
                            & \lesssim \int_{2^m}^\infty \frac{\capac(\Omega^c \cap \overline B_t ,B_{4t})}{t^{n + \gamma - 2}} \frac{\mathrm dt}{t}    \\[6pt]
                            & = \sum_{k = m + 1}^\infty \int_{2^{k - 1}}^{2^k} \frac{\capac(\Omega^c \cap \overline B_t, B_{4t})}{t^{n + \gamma - 2}} \frac{\mathrm dt}{t}   \\[6pt]
                            & \leq \sum_{k = m + 1}^\infty \int_{2^{k - 1}}^{2^k} \frac{\capac(\Omega^c \cap \overline B_t, B_{2^{k + 1}})}{t^{n + \gamma - 2}} \frac{\mathrm dt}{t}   \\[6pt]
                            & \leq \sum_{k = m + 1}^\infty \capac(\Omega^c \cap \overline B_{2^k},B_{2^{k + 1}}) \int_{2^{k - 1}}^{2^k} \frac{1}{t^{n + \gamma - 2}} \frac{\mathrm dt}{t}   \\[6pt]
                            &\lesssim \sum_{k = m + 1}^\infty (2^{-k})^{n + \gamma - 2} \capac(\Omega^c \cap \overline B_{2^k},B_{2^{k + 1}}) \eqqcolon T. 
        \end{align}
    Defining $K \coloneqq (2^{-(m + 1)})^{n + \gamma - 2} \capac(\Omega^c \cap B_{2^(m + 1)},B_{2^{m + 2}})$ and using that $B_{2^k} \subset A_k \cup B_{2^{k - 1}}$ we obtain
        \begin{align}
            T   & = K + \sum_{k = m + 2}^\infty (2^{-k})^{n + \gamma - 2} \capac(\Omega^c \cap \overline B_{2^k},B_{2^{k + 1}})\\[6pt]
                & \leq K + \sum_{k = m + 2}^\infty (2^{-k})^{n + \gamma - 2} \capac(\Omega^c \cap A_k,B_{2^{k + 1}})    \\[6pt]
                & \quad + \sum_{k = m + 2}^\infty (2^{-k})^{n + \gamma - 2} \capac(\Omega^c \cap \overline B_{2^{k - 1}},B_{2^{k + 1}}) \\[6pt]
                & = K + \WW^\Sigma(m + 2,\Omega) + \sum_{k = m + 1}^\infty (2^{-(k + 1)})^{n + \gamma - 2} \capac(\Omega^c \cap \overline B_{2^k},B_{2^{k + 2}}) \\[6pt]
                & = K + \WW^\Sigma(m + 2,\Omega) + 2^{2 - (n + \gamma)} \sum_{k = m + 1}^\infty (2^{-k})^{n + \gamma - 2} \capac(\Omega^c \cap \overline B_{2^k},B_{2^{k + 2}})  \\[6pt]
                & \leq K + \WW^\Sigma(m + 2,\Omega) + 2^{2 - (n + \gamma)} \sum_{k = m + 1}^\infty (2^{-k})^{n + \gamma - 2} \capac(\Omega^c \cap \overline B_{2^k},B_{2^{k + 1}})  \\[6pt]
                & = K + 2^{2 - (n + \gamma)} T + \WW^\Sigma(m + 2,\Omega).
        \end{align}
    After absorption, we therefore find that
        \begin{align}
            \WW(2^m,\Omega) \lesssim T \lesssim 1 + \WW^\Sigma(m + 2,\Omega),
        \end{align}
    and therefore the claim. This establishes the equivalence \textit{(4)} $\Longleftrightarrow$ \textit{(5)} in Theorem \ref{maintheorem}.
    
    \medskip

    We end this section by recalling a bound on potentials in terms of the Wiener integral. The proof can be found in \cite[Lemma 6.25]{HeKiMa:2006} (see also \cite[Corollary 1]{Mazja:1976})

     \begin{lemma} \label{lemma:exponential_estimate}
        Let $x_0 \in \partial \Omega$ and let $0 < r \leq R$. Then there exists $c = c(n,\gamma) > 0$ such that
            \begin{align}
                1 - \hat \R^1[\Omega^c \cap \overline B_R(x_0),B_{2R}(x_0)](x) \leq \exp\Big( -c \int_r^R \frac{\capac(\Omega^c \cap B_t(x_0), B_{2t}(x_0))}{\capac(B_t(x_0), B_{2t}(x_0))} \frac{\mathrm dt}{t} \,\Big)
            \end{align}
        for all $x \in B_{\frac r2}(x_0)$.
    \end{lemma}

\section{Proof of the Main Theorem}

    We separate the proof of the main theorem into multiple propositions. The first proposition, similar to the case of uniformly elliptic operators (\cite[Lemma 2.1] {Abdulla:2007}, \cite[Theorem 1.1, (1)$\Longleftrightarrow$(2)] {Abdulla:2012}), connects the well-posedness of \eqref{eq:dirichlet_problem} and the $\A$-harmonic measure of $\infty$.

    \begin{prop} \label{prop:harmonic_measure::uniqueness}
        The following assertions are equivalent.
            \begin{enumerate}
                \item $\infty$ is regular (resp. irregular) for $\Omega$.
                \item For every bounded Borel measurable function $f \colon \partial \Omega \longrightarrow \IR$ there exists a unique generalised solution of \eqref{eq:dirichlet_problem} (resp. there exist infinitely many generalised solutions of \eqref{eq:dirichlet_problem}).
            \end{enumerate}
    \end{prop}
        \proofb Assume that $\infty$ is regular, i.e. $\H^\Omega_{\mathbf 1_\infty} = 0$. Let $f \colon \partial \Omega \longrightarrow \IR$ be a bounded Borel measurable function and let $u, \tilde u \colon \Omega \longrightarrow \IR$ be generalised solutions of \eqref{eq:dirichlet_problem}. Then $v \coloneqq u - \tilde u$ is a generalised solution of \eqref{eq:dirichlet_problem} with trivial boundary data and thus there exists $M \in \IR$ such that $v = M \H^\Omega_{\mathbf 1_\infty} = 0$ by assumption. Therefore $u = \tilde u$.

        \medskip

        We now assume that $\infty$ is irregular, i.e. $\H_{\mathbf 1_\infty} > 0$ in $\Omega$. Let $f \colon \partial \Omega \longrightarrow \IR$ be a bounded Borel measurable function and let $u \colon \Omega \longrightarrow \IR$ be a generalised solution of \eqref{eq:dirichlet_problem}. For every $M \in \IR$, the function $v_M \coloneqq u + M \H_{\mathbf 1_\infty}$ is also a generalised solution to \eqref{eq:dirichlet_problem}. Since $v_M \neq v_{\tilde M}$ for every $M \neq \tilde M$, the claim follows. \proofe

    The next proposition draws a connection between the well-posedness of \eqref{eq:dirichlet_problem} and the condition \eqref{eq:defi_bc}
(\cite[Lemma 2.1] {Abdulla:2007}, \cite[Theorem 1.1, (2)$\Leftrightarrow$ (3)] {Abdulla:2012}).

    \begin{prop} \label{prop:uniqueness::regular_point}
        The following assertions are equivalent.
            \begin{enumerate}
                \item For every bounded Borel measurable function $f \colon \partial \Omega \longrightarrow \IR$ there exists a unique generalised solution of \eqref{eq:dirichlet_problem} (resp. there exist infinitely many generalised solutions of \eqref{eq:dirichlet_problem}).
                \item \eqref{eq:defi_bc} is satisfied (resp. is not satisfied).
            \end{enumerate}
    \end{prop}
        \proofb Assume that there exists a unique generalised solution of \eqref{eq:dirichlet_problem} for every bounded Borel measurable function $f \colon \partial \Omega \longrightarrow \IR$. We want to show that \eqref{eq:defi_bc} is satisfied. Assume, on the contrary, that \eqref{eq:defi_bc} is not satisfied. Thus, we would find a bounded Borel measurable function $f \colon \partial \Omega \longrightarrow \IR$ and a generalised solution $u \colon \Omega \longrightarrow \IR$ of \eqref{eq:dirichlet_problem} such that either $f^* < u^*$ or $f_* > u_*$. However, using Lemma \ref{lem:construction_generalised_solution_good_limit}, there exists a generalised solution $v \colon \Omega \longrightarrow \IR$ of \eqref{eq:dirichlet_problem} such that $f_* \leq v_* \leq v^* \leq f^*$, which contradicts uniqueness.

        \medskip

        Assume now that \eqref{eq:defi_bc} is not satisfied. Thus, there exists a generalised solution $u \colon \Omega \longrightarrow \IR$ and a bounded Borel measurable function $f \colon \partial \Omega \longrightarrow \IR$ such that either $f^* < u^*$ or $f_* > u_*$. Now choose $\overline f \in \IR$ such that $f_* \leq \overline f \leq f^*$. Then the generalised solution $v \colon \Omega \longrightarrow \IR$ of \eqref{eq:dirichlet_problem} with boundary data $f$ extended via $f(\infty) =\overline f$ satisfies $u \neq v \coloneqq \H^\Omega_f$ (see Lemma \ref{lem:construction_generalised_solution_good_limit}), from which we deduce that $w \coloneqq u - v$ is a non-trivial generalised solution of \eqref{eq:dirichlet_problem} with trivial boundary data. Therefore $\H^\Omega_{\mathbf 1_\infty}$ is non trivial, and thus there exist infinitely many solutions of \eqref{eq:dirichlet_problem}, namely $u_M \coloneqq u + M \H^\Omega_{\mathbf 1_\infty}$ for every $M \in \IR$. \proofe

    The following proposition connects the Wiener integral and the $\A$-harmonic measure of $\infty$. Implementing the method developed in \cite[Theorem 1.1] {Abdulla:2007} and \cite[Theorem 1.1, (1)$\Rightarrow$ (4)] {Abdulla:2012}, it is demonstrated that the convergence of the Wiener integral allows for suitable lower bounds of potentials (in terms of their respective capacities) of sets that exhaust the exterior set $\Omega^c$. These potentials are then used to bound the $\A$-harmonic measure of $\infty$ from below.

    \begin{prop} \label{prop:Wiener_converges::irregular}
       Let $\Omega \subset \IRn$ be open such that there exists an unbounded connected component. Assume that the Wiener integral converges. Then $\infty$ is irregular for $\Omega$.
    \end{prop}
        \proofb We show that $(\H^\Omega_{\mathbf 1_\infty})^* = 1$. Since the PWB-method is order-preserving, we have $(\H^\Omega_{\mathbf 1_\infty})^* \leq 1$, so we aim to establish the reverse inequality. Using the results from Section \ref{subsec:wiener_integral}, we know that the Wiener sum $\WW^\Sigma$ is finite (without loss of generality, we take the centre point to be the origin). Denote $B_t \coloneqq B_t(0)$ and let $\varepsilon > 0$. Since the Wiener sum converges, we find $m \in \IN$ such that
            \begin{align} \label{eq:definition_exterior_annulli}
                \sum_{k = m + 1}^\infty (2^{-k})^{n + \gamma - 2} \capac(E_k, B_{2^{k + 1}}) < \varepsilon,    && E_k \coloneqq \Omega^c \cap \{2^{k - 1} \leq |x| \leq 2^k\} \subset \overline B_{2^k}.
            \end{align}
        Since $\Omega$ has at least one unbounded connected component $\Omega_e \subseteq \Omega$, we can choose $m$ sufficiently large such that there exists $x^* \in \Omega_e \cap B_{2^{m - 4}}$. We note that once we have found $m$ sufficiently large and $x^* \in \Omega \cap B_{2^{m - 4}}$, the point $x^*$ is independent of $m$ in that for larger $m$ we can still choose $x^*$ as before. In particular, constants that depend on $|x^*|$ do not depend on $m$.
        
        \medskip

        The proof is decomposed into multiple steps. In the first step we construct an $\A$-superharmonic function $\vartheta \colon \IRn \longrightarrow \IR$ such that the liminf $\vartheta_* \lesssim \varepsilon$ at $\infty$ and $\vartheta \geq 1$ quasi-everywhere on $\partial \Omega^m$, where in $\Omega^m \coloneqq \big(\bigcup_{k = m + 1}^\infty E_k \big)^c$. We then modify $\vartheta \rightsquigarrow \varpi$ in such a way that the liminf $\varpi_* \lesssim \varepsilon$ and $\varpi \geq 1$ everywhere on $\partial \Omega^m$. In fact $\varpi \geq 1$ in $(\Omega^m)^c$. Finally, having obtained these properties, the function $\varpi + \hat \R^1[B_{2^m},\IRn]$ is in the upper class of the function $1 - \mathbf 1_\infty$, which in turn is used to bound the harmonic measure.

        \medskip
        
        \textit{Step 1: Construction of potential.} Given $M \in \IN$ such that $M > m$, we define $\vartheta_M \colon \IRn \longrightarrow \IR$ via
            \begin{align}
                \vartheta_M \coloneqq \sum_{k = m + 1}^M \hat \R^1[E_k,\IRn].
            \end{align}
        Then $\vartheta_M$ is $\A$-harmonic in $\Omega^m$ for all $M > m$. Moreover, $\{\vartheta_M\}_M$ is an increasing sequence of non-negative $\A$-harmonic functions in $\Omega^m$. We claim that the pointwise limit function $\vartheta \colon \IRn \longrightarrow \IR$ is an $\A$-harmonic function in $\Omega^m$. We note that
            \begin{align}
                \vartheta(x) = \lim_{M \to \infty} \sum_{k = m + 1}^M \hat \R^1[E_k,\IRn](x) = \sup_{M > m} \vartheta_M(x)  && \text{for all } x \in \Omega^m.
            \end{align}
        Let $x^* \in \Omega \cap B_{2^{m - 4}}$. Let $M - 2> k > m$, $G_M(x^*, \filll) \in L^1(B_{2^M}, \mu_\gamma)$ be the Greens function (with singularity in $x^*$) and let $\nu_{E_k}^M$ be the capacitary distribution of $\hat \R^1[E_k, B_{2^M}]$. We use the Greens representation formula \eqref{eq:greens_representation_formula}, \eqref{eq:greens_function_estimate_local}, Lemma \ref{lemma:capacity_properties} and \eqref{eq:structural_capacity_estimate} to obtain
            \begin{align}
\label{eq:wiener_potential1}    \hat \R^1[E_k,\IRn](x^*)    & = \lim_{M \to \infty} \int_{E_k} G_M(x^*,y) \dd \nu_{E_k}^M(y)    \\[6pt]
                        & \lesssim \liminf_{M \to \infty} \int_{E_k} \int_{|x^* - y|}^{2^M} \frac{1}{\capac(B_t(x^*), B_{2t}(x^*))} \frac{\mathrm dt}{t} \dd \nu_{E_k}^M(y) \\[6pt]
                        & \leq \liminf_{M \to \infty} \int_{E_k} \int_{|x^* - y|}^\infty \frac{1}{\capac(B_t(x^*), B_{2t}(x^*))} \frac{\mathrm dt}{t} \dd \nu_{E_k}^M(y)  \\[6pt]
                        & \lesssim  \capac(E_k, B_{2^{k + 1}}) \int_{2^{k - 1} - |x^*|}^\infty \frac{1}{{\capac(B_t(x^*), B_{2t}(x^*))}} \frac{\mathrm dt}{t} \\[6pt]
                        & \lesssim \capac(E_k, B_{2^{k + 1}}) \int_{2^{k - 1} - |x^*|}^\infty \frac{t}{{\mu_\gamma(B_t(x^*))}} \dd t \\[6pt]
                        & \lesssim \capac(E_k, B_{2^{k + 1}}) \int_{2^{k - 1} - |x^*|}^\infty \frac{t}{(t - |x^*|)^{n + \gamma}} \dd t    \\[6pt]
                        & = \capac(E_k, B_{2^{k + 1}}) \bigg[\frac{2^{2 - (n + \gamma)}}{n + \gamma - 2} (2^{k - 2} - |x^*|)^{2 - (n + \gamma)}   \\[6pt]
                        & \qquad + |x^*| \, \frac{2^{1 - (n + \gamma)}}{n + \gamma - 1} (2^{k - 2} - |x^*|)^{1 - (n + \gamma)} \bigg] \\[6pt]
\label{eq:wiener_potential2} & \lesssim \capac(E_k, B_{2^{k + 1}}) \, 2^{-k (n + \gamma - 2)},
            \end{align}
        where in \eqref{eq:wiener_potential2} it was used that $|x^*| \leq 2^{m - 4}$. Using that the Wiener integral converges, we therefore obtain $\vartheta(x^*) \lesssim \varepsilon$, and hence $\vartheta$ is $\A$-superharmonic in $\IRn$ and $\A$-harmonic in $\Omega$. In particular, we find that the liminf at $\infty$ is estimated via $\vartheta_* \lesssim \varepsilon$ using the comparison principle. By construction, we find that $\vartheta \geq 1$ quasi everywhere on $\partial \Omega^m$ as desired.
        
        \medskip
        
        \textit{Step 2: Modification of $\vartheta$.} We aim to modify $\vartheta \rightsquigarrow \varpi$ in such a way that $\varpi_* \lesssim \varepsilon$ at $\infty$ and $\varpi \geq 1$ everywhere on $\partial \Omega^m$. To do so, we define
            \begin{align}
                \tilde E_k \coloneqq \{x \in E_k : \vartheta(x) < 1\}  && \text{for all } k \in \IN, \ k > m.
            \end{align}
        Since the set of irregular boundary points is a set of capacity $0$, we find an open set $\tilde E_k \subset U_k \subset B_{2^{k + 1}}$ such that $\capac(\overline U_k, B_{2^{k + 1}}) < \varepsilon$. By definition of capacity, we can, without loss of generality, assume that $U_k \subset B_{2^{k - \frac 32}}^c$. We now consider the $\A$-superharmonic function $\varsigma \colon \IRn \longrightarrow \IR$ defined via
            \begin{align}
                \varsigma \coloneqq \sum_{k = m + 1}^\infty \hat \R^1[\overline U_k, \IRn]
            \end{align}
        We repeat the computations from \eqref{eq:wiener_potential1} -- \eqref{eq:wiener_potential2} to obtain
            \begin{align}
                \hat \R^1[\overline U_k,\IRn](x^*)  & \lesssim \liminf_{M \to \infty} \int_{\partial U_k} \int_{|x^* - y|}^\infty \frac{1}{\capac(B_t(x^*), B_{2t}(x^*))} \frac{\mathrm dt}{t} \dd \nu_{E_k}^M(y)  \\[6pt]
                                                    & \lesssim  \capac(\overline U_k, B_{2^{k + 1}}) \int_{2^{k - \frac 32} - |x^*|}^\infty \frac{1}{{\capac(B_t(x^*), B_{2t}(x^*))}} \frac{\mathrm dt}{t} \\[6pt]
                                                    & \lesssim \varepsilon \int_{2^{k - \frac32} - |x^*|}^\infty \frac{t}{{\mu_\gamma(B_t(x^*))}} \dd t \\[6pt]
                                                    & \lesssim \varepsilon \int_{2^{k - \frac 32} - |x^*|}^\infty \frac{t}{(t - |x^*|)^{n + \gamma}} \dd t    \\[6pt]
                                                    & = \varepsilon \big[2^{2 - (n + \gamma)} (2^{k - \frac 52} - |x^*|)^{2 - (n + \gamma)} + |x^*| 2^{1 - (n + \gamma)} (2^{k - \frac 52} - |x^*|)^{1 - (n + \gamma)} \big] \\[6pt]
\label{eq:exponential_estimate2}                    & \lesssim \varepsilon  2^{-k (n + \gamma - 2)}.
            \end{align}
        In total, we find that
            \begin{align}
                \varsigma(x^*) \lesssim \varepsilon \sum_{k = m + 1}^\infty (2^{-k})^{n + \gamma - 2} \lesssim \varepsilon.
            \end{align}
        Now by definition, the function $\varpi \coloneqq \vartheta + \varsigma \colon \IRn \longrightarrow \IR$ is $\A$-superharmonic in $\IRn$, satisfying $\varpi \geq 1$ everywhere on $(\Omega^m)^c$, in particular on $\partial \Omega^m$, and $\varpi(x_*)\lesssim\varepsilon$. Since $\varpi$ is $\A$-superharmonic in $\Omega^m$, by applying $\A$-superharmonic minimum principle in $\Omega^m\cap B_R$ we conclude that for all sufficiently large $R$, there must be a point $x_R\in \Omega_e\cap B_R$ such that $\varpi(x_R) \lesssim \varepsilon$. Hence, we deduce that $\varpi_* \lesssim \varepsilon$ at $\infty$. 
        
        \medskip
        
        \textit{Step 3: Conclusion.} Using step 2, the function $\varpi + \hat \R^1[B_{2^m}, \IRn]$ is now in the upper class of the boundary function $1 - \mathbf{1}_{\infty}$. Thus, we find that
            \begin{align}
                0 \leq 1 - \H^\Omega_{\mathbf 1_\infty} = \H^{\Omega}_{1 - \mathbf 1_\infty} \leq \varpi + \hat \R^1[B_{2^m},\IRn]  && \text{in } \Omega,
            \end{align}
        from which we obtain
            \begin{align}
                0 \leq \liminf_{\substack{x \to \infty \\ x \in \Omega}} 1 - \H^\Omega_{\mathbf 1_\infty}(x) \leq \varpi_* + 0 \lesssim \varepsilon.
            \end{align}
        Thus, there exists a constant $C = C(n, \gamma) > 0$ such that
            \begin{align}
                1 - C\varepsilon \leq \limsup_{\substack{x \to \infty \\ x \in \Omega}} \H^\Omega_{\mathbf 1_\infty}.
            \end{align}
        Since $\varepsilon > 0$ was arbitrary, the claim follows. \proofe

    The following proposition shows the other direction of Proposition \ref{prop:Wiener_converges::irregular}. It is based on local bounds in terms of the Wiener integral that are obtained via a well-known iteration process of potentials (see \cite[Lemma 6.25]{HeKiMa:2006}, see also \cite[Corollary 1]{Mazja:1976}). In contrast to the preceding proposition, we do not require $\Omega$ to have an unbounded connected component.

    \begin{prop} \label{prop:Wiener_diverges::regular}
      Assume that the Wiener integral diverges. Then $\infty$ is regular for $\Omega$.  
    \end{prop}
        \proofb Let $x_0 \in \partial \Omega$ and let $0 < r < R$. We write $B_R \coloneqq B_R(x_0)$. We define the function $v_R \colon \IRn \longrightarrow \IR$ via
            \begin{align}
                v_R \coloneqq \left\{ \begin{array}{ll}
                     1 - \hat \R^1[\Omega^c \cap \overline B_R,B_{2R}]  & \quad \text{in } B_{2R},  \\[6pt]
                     1                                                  & \quad \text{in } B_{2R}^c.
                \end{array}\right.
            \end{align}
        Using the pasting lemma for $\A$-superharmonic functions, the function $v_R$ is $\A$-superharmonic in $\Omega$ and therefore in the upper class of the boundary function $\mathbf 1_\infty$ for every $R > 0$. Using Lemma \ref{lemma:exponential_estimate}, there exists $c = c(n,\gamma) > 0$ such that
            \begin{align}
               0 \leq \H_{\mathbf 1_\infty}^\Omega(x) \leq v_R(x) \leq \exp\Big( -c \int_r^R \frac{\capac(\Omega^c \cap B_t, B_{2t})}{\capac(B_t, B_{2t})} \frac{\mathrm dt}{t} \,\Big)
            \end{align}
        for all $x \in B_{\frac r2} \cap \Omega$. Using that the Wiener integral diverges, in the limit $R \to \infty$ we find that $\H^\Omega_{\mathbf 1_\infty} = 0$ in $\Omega \cap B_{\frac r2}$. Since $r > 0$ was arbitrary, the claim follows. \proofe

    The next proposition establishes the last missing link \textit{(1)} $\Longleftrightarrow$ \textit{(3)} in Theorem \ref{maintheorem}. It connects the irregularity of $\infty$ to $\A$-thinness at $\infty$ of the exterior set (see \cite[Theorem 1.1] {Abdulla:2012} for the case of uniformly elliptic operators).

    \begin{prop} The following assertions are equivalent.
       \begin{enumerate}
           \item $\Omega^c$ is $\A$-thick (resp. $\A$-thin) at $\infty$.
           \item $\infty$ is regular (resp. irregular) for $\Omega$.
       \end{enumerate}
    \end{prop}
    
    The proof of this proposition is identical to the proof given in the case of uniformly elliptic operators (see \cite[Theorem 1.1, (1) $\Longleftrightarrow$ (5)] {Abdulla:2012}). 

\bibliographystyle{alpha}
\bibliography{bibliography}

\end{document}

%% file: header.tex
\usepackage[utf8]{inputenc}
\usepackage[english]{babel}
\usepackage{amsfonts}
\usepackage{amssymb}
\usepackage{amsmath}
	\usepackage{mathrsfs}						%mathematische Sonderzeichen
	\usepackage{mathtools}
	\usepackage{bbm}
	\usepackage[]{units} 						%Bessere Darstellung im mathmode
\usepackage{amsthm}								%'newtheorem' Umgebung
\usepackage{graphicx}							%Bilder einbinden
\usepackage{caption}
\usepackage{pgf,tikz}							%Tikz
	\usepackage{mathrsfs}
\usepackage{color}								%Farben
\usepackage{soul} 								% strike \st{}, underline \ul{}, highlight \hl{}
\usepackage{paralist}							%Aufzählung ändern
\usepackage{setspace}							%
\usepackage{esint}								%Mittelwertintegral
\usepackage{microtype}							%Mikrotypografie
\usepackage{framed}								%Hintergründe
\usepackage[pdfencoding=auto, psdextra, unicode]{hyperref}%Referenzieren
\usepackage{enumitem}
\usepackage{geometry}							%Satzspiegel
plus 6pt minus 12pt
\allowdisplaybreaks                         %Align mit Seitenumbrüchen

\hypersetup{pdfborder = {0 0 0}}			%Unsichtbarer Querverweis	

\numberwithin{equation}{section}			%Nummerierung der Gleichungen
\numberwithin{figure}{section}				%Nummerierung der Abbildungen

\definecolor{shadecolor}{RGB}{238,233,233}	%Hintergrundfarbe

\renewcommand{\div}{\textrm{div}}			%Divergenz
\newtheoremstyle{script}					%Name
  {}										%Abstand nach oben
  {}										%Abstand nach unten
  {\itshape}								%Schrift
  {}										%Theoremkopf Einrückung
  {\sffamily\bfseries}						%Theoremkopf Schrift
  {}										%Punktuation nach Theoremkopf
  {5pt}										%Abstand nach Theoremkopf
  {}										%Theoremkopfspezfikation

\theoremstyle{script}
\newtheorem{Thrm}{Theorem}[section]
\newtheorem{Defi}[Thrm]{Definition}

\newtheorem{Coro}[Thrm]{Corollary}
\newtheorem{Lemma}[Thrm]{Lemma}
\newtheorem{Prop}[Thrm]{Proposition}

\DeclareMathOperator*{\Essliminf}{ess\, lim \, inf}					%Essentieller Limes Inferior

\newenvironment{defi}
	{\begin{Defi}}
	{\end{Defi}}
	
\newenvironment{thrm}
	{\begin{Thrm}}
	{\end{Thrm}}

\newenvironment{lemma}
	{\begin{Lemma}}
	{\end{Lemma}}

\newenvironment{prop}
	{\begin{Prop}}		
	{\end{Prop}}

\newcommand{\IN}{\mathbbm{N}}

\newcommand{\IR}{\mathbbm{R}}

\newcommand{\A}{\mathcal{A}}

\newcommand{\G}{\mathcal{G}}	
\renewcommand{\H}{\mathcal{H}}

\newcommand{\R}{\mathcal{R}}

\newcommand{\WW}{\mathfrak{W}}

\newcommand{\IRn}{{\IR^{\kern-.7pt n}}}							%Euklidische Räume
\newcommand{\IRN}{{\IR^{\kern-.7pt N}}}

\newcommand{\proofb}{\noindent\textit{Proof:}\ }
\newcommand{\proofe}{\hfill $\blacksquare$}

\newcommand{\filll}{\ \cdot\ }									%Operator - Platzhalter

\newcommand{\dd}{\, \textrm{d}}									%Differential
\newcommand{\supp}{\textrm{supp}}								%Träger
\newcommand{\loc}{{\textrm{loc}}}

\newcommand{\capac}{\operatorname{cap}}